\documentclass[12pt]{article}
\usepackage{amsmath,amssymb}

\def\seq#1#2#3{#1_{#2},\,\ldots,#1_{#3}}

\def\VV{{\underline{V}}}
\def\vv{{\underline{v}}}

\def\tt{{\underline{t}}}
\def\TT{{\underline{T}}}

\def\mm{\underline{m}}

\def\uu{\underline{u}}
\def\ll{\underline{\ell}}

\def\1{\underline{1}}

\def\P{\Bbb P}

\def\Z{\Bbb Z}

\def\C{\Bbb C}

\def\CP{\Bbb C\Bbb P}
\def\OO{{\cal O}}

\newtheorem{theorem}{Theorem}

\newenvironment{definition}
{\smallskip\noindent{\bf Definition\/}:}{\smallskip\par}

\newenvironment{examples}
{\smallskip\noindent{\bf Examples\/}.}{\medskip\par}
\newenvironment{remark}
{\smallskip\noindent{\bf Remark\/}.}{\smallskip\par}
\newenvironment{remarks}
{\smallskip\noindent{\bf Remarks\/}.}{\smallskip\par}
\newenvironment{proof}
{\noindent{\bf Proof\/}.}{{ $\Box$}\smallskip\par}

\title{On Poincar\'e series of filtrations on equivariant functions
of two variables
\footnote{Math. Subject Class. 14B05, 16W70, 16W22.}
}

\author{
A.~Campillo
\and F.~Delgado \and S.M.~Gusein-Zade
\thanks{Partially supported by the grant MCYT: MTM2004-00958.
The third author
partially supported by the grants RFBR--04--01--00762, NSh--4719.2006.1.}
}

\date{}
\begin{document}
\def\eps{\varepsilon}

\maketitle

\begin{abstract}
Let a finite group $G$ act on the complex plane $(\C^2, 0)$. We consider
multi-index filtrations on the spaces of germs of holomorphic functions
of two variables equivariant with respect to 1-dimensional representations
of the group $G$ defined by components of a modification of the complex
plane $\C^2$ at the origin or by branches of a $G$-invariant plane curve
singularity $(C,0)\subset(\C^2,0)$. We give formulae for the Poincar\'e series
of these filtrations. In particular, this gives a new method to obtain the Poincar\'e series of analogous filtrations on the rings of germs of functions on quotient surface singularities.
\end{abstract}

\section*{Introduction}\label{sec0}
In  \cite{CDGa}, \cite{DG} (see also \cite{CDGb} for another proof
of the result of \cite{CDGa}) there were computed the Poincar\'e
series (in several variables) of the so called divisorial (multi-index)
filtrations on the ring ${\cal O}_{\C^2,0}$ of germs of holomorphic functions on
the complex plane $\C^2$ at the origin and of the filtration
defined by orders of a function on branches of a germ of a plane
curve $(C,0)\subset(\C^2,0)$. In particular it was shown that the
Poincar\'e series of the latter filtration coincides with the
(multi-variable) Alexander polynomial of the (algebraic) link
$(S_\eps^3, C\cap S_\eps^3)$ corresponding to the curve $(C,0)$.
Here we give versions of these results
for the equivariant case, i.e. when there is an action (a
representation) of a finite group $G$ on the complex plane $\C^2$ and the
corresponding filtrations are considered on subspaces of germs of
functions equivariant with respect to 1-dimensional
representations $G\to \C^*=\mbox{\bf GL}(\C, 1)$ of the group $G$. There is a
problem
to define an equivariant version of the monodromy zeta function. One reason to
study "equivariant" Poincare series of filtrations is the hope that they may
give a hint for such a definition (at least for curve singularities).

For a finite group $G$, the set $G^*$ of 1-dimentional representations
$G\to\C^*=\mbox{\bf GL}(\C,1)$ of the group $G$ is a group itself.
Let $R(G)$ be the ring of (complex virtual) representations of
the group $G$ and let $R_1(G)$ be its subring generated by
1-dimensional representations. Elements of the ring $R_1(G)$
are formal sums $\sum n_\alpha \alpha$ over all different
1-dimensional representations $\alpha$ of the group $G$, $n_\alpha\in \Z$. The
multiplication is defined by the tensor product. The rings $R(G)$ and $R_1(G)$
coincide iff the group $G$ is abelian.

Let $(V,0)$ be a germ of a complex analytic variety
with an action of a finite group $G$. The group $G$ acts
on the ring ${\cal O}_{V,0}$ of germs of functions on $(V,0)$ by
the formula $g^*f(x)=f(g^{-1}(x))$ ($f\in {\cal O}_{V,0}$, $g\in G$,
$x\in(V,0)$). Let $\alpha:G\to \C^*$ be a 1-dimensional representation of the group $G$

\begin{definition}
A germ $f\in {\cal O}_{V,0}$ is {\em equivariant} with respect to
the representation $\alpha$ if $g^*f=\alpha(g)\cdot f$ for $g\in
G$.
\end{definition}

The set of germs of functions on $(V,0)$ equivariant with respect to the
representation $\alpha$ is a vector subspace
${\cal O}^\alpha_{V,0} \subset {\cal O}_{V,0}$. If $\alpha_1\ne\alpha_2$,
one has ${\cal O}^{\alpha_1}_{V,0}\cap{\cal O}^{\alpha_2}_{V,0} =
\{0\}$. Product of functions equivariant with respect to
representations $\alpha_1$ and $\alpha_2$ is a function
equivariant with respect to the representation
$\alpha_1\cdot\alpha_2$.
If the group $G$ is abelian,
${\cal O}_{V,0} = \bigoplus\limits_{\alpha}{\cal O}^\alpha_{V,0}$.

For a vector space $A$ (finite or infinite dimensional, e.g., for
${\cal O}_{V,0}$ or ${\cal O}^\alpha_{V,0}$), a $1$-index (decreasing)
filtration on $A$ is defined by a function $v: A \to \Z_{\ge 0}
\cup \{\infty\}$, such that $v(f_1+f_2)\ge \min\{v(f_1),
v(f_2)\}$, $f_i\in A$, $i=1,2$. An $s$-index filtration on $A$ is
defined by $s$ such functions $v_i$, $i=1, \ldots, s$. For
$\vv=(v_1, \ldots, v_s)\in \Z^s$, the corresponding subspace
$J(\vv)$ is defined as $\{a\in A: v_i(a)\ge v_i, i=1,\ldots, s\}$.
(To describe the filtration it is sufficient to define the
subspaces $J(\vv)$ only for $\vv\in\Z^s_{\ge 0}$ (i.e., $v_i\ge 0$
for $i=1, \ldots, s$), however it is convenient to suppose
$J(\vv)$ to be defined for all $\vv$ from the lattice $\Z^s$.)

Suppose that all the factor
spaces $J(\vv)/J(\vv+\1)$ ($\1=(1,\ldots,1)$, $\vv\in\Z^s$) are
finite dimensional. (This is equivalent to say that $\dim
A/J(\vv)$ is finite for any $\vv\in\Z^s_{\ge 0}$.) Let $d(\vv) :=
\dim J(\vv)/J(\vv+\1)$, $\vv\in \Z^r$. The {\em Poincar\'e series}
of the multi-index filtration $\{J(\vv)\}$ is the series
$P(t_1,\ldots, t_s)\in\Z[[t_1,\ldots, t_s]]$ defined by the
formula
\begin{equation}\label{eq1}
P(\tt)=\frac{\left(\sum \limits_{\vv\in\Z^s}
d(\vv)\cdot\tt^{\vv}\right)\prod\limits_{i=1}^s(t_i-1)}
{\tt^{\1}-1}\,,
\end{equation}
where $\tt^{\vv} = t_1^{v_1}\cdot \ldots \cdot t_s^{v_s}$
($\tt=(t_1,\ldots,t_s)$, $\vv=(v_1,\ldots,v_s)$); see, e.g.,
\cite{CDK}, \cite{CDGa}.

\begin{remark}
Pay attention that the sum in the numerator is over all
$\vv\in\Z^s$, not only over those from $\Z^s_{\ge 0}$. This sum
contains monomials with negative powers of variables (for $s\ge
2$), however, being multiplied by $\prod\limits_{i=1}^s(t_i-1)$,
it becomes a power series, i.e. an element of $\Z[[t_1,\ldots,
t_s]]$.
\end{remark}

If there is a linear action
(a representation) of a group $G$ on the space $A$ and the
filtration $\{J(\vv)\}$ ($\vv\in\Z^s$) is invariant with respect
to this representation (i.e. $g^*J(\vv)=J(\vv)$ for $g\in G$), one can define
the equivariant Poincar\'e series of the filtration $\{J(\vv)\}$
as an element of $R(G)[[t_1,\ldots, t_s]]$. For that $d(v)$ in
(\ref{eq1}) should be substitute by $[J(\vv)/J(\vv+\1)]$ considered
as a $G$-module (i.e. an element of the ring $R(G)$). Filtrations on the
ring ${\cal O}_{V,0}$ considered here are, generally speaking,
not $G$-invariant: the action of the group permutes valuations.
This difficulty can be (partially) avoided by considering
filtrations on the spaces of equivariant functions.
This way an equivariant Poincar\'e series can be defined as an
element of $R_1(G)[[t_1,\ldots, t_s]]$.

In the case when $A$ is the ring ${\cal O}_{V,0}$ of germs of
functions on $(V,0)$ or a subspace of it, one can define a notion
of integration with respect to the Euler characteristic over the
projectivization $\P A$ of the space $A$ (see, e.g., \cite{CDGb},
\cite{CDGc}; there the notion is defined for $A={\cal O}_{V,0}$,
however there is no essential difference with the general case
$A\subset {\cal O}_{V,0}$: the role of the jet spaces $J^k_{V,0}=
{\cal O}_{V,0}/{\frak m}^{k+1}$ ($\frak m$ is the maximal ideal
in the ring ${\cal O}_{V,0}$) is played by the spaces $JA^k_{V,0}=
A/{\frak m}^{k+1}\cap A$). This notion generalizes the usual notion
of integration with respect to the Euler characteristic (\cite{Viro})
and was inspired by the notion of motivic integration (see, e.g.,
\cite{DL}).

\begin{definition}
The (multi-index) filtration $\{J(\vv)\}$ on a subspace
$A\subset {\cal O}_{V,0}$ is {\em finitely determined} if for each
$\vv\in\Z^s$ there exists $N\ge 0$ such that
${\frak m}^N\cap A \subset J(\vv)$.
\end{definition}

If the filtration $\{J(\vv)\}$ is finitely determined, one can show that
$$
P(\tt)=\int\limits_{\P A} \tt^{\,\vv(f)}d\chi\,;
$$
see, e.g., \cite{CDGb}, \cite{CDGc}.

Here we compute the Poincar\'e series for some multi-index filtrations
(divisorial ones and those defined by invariant curves) on the space of
equivariant function-germs on the plane $(\C^2, 0)$ with
a finite group action. For the trivial (1-dimensional) representation of
the group, this gives formulae for divisorial filtrations on the quotient
surface singularity $(\C^2/G, 0)$ and for the filtrations defined by curves
on them somewhat different from those in $\cite{CDGc}$, $\cite{CDGd}$.

\section{Divisorial filtrations on the spaces of equivariant
functions}\label{sec1}
Let $\pi:(\widetilde V, D)\to (V,0)$ be a resolution of the germ $(V,0)$.
Here $\widetilde V$ is a smooth manifold, $\pi$ is a proper analytic map
which is an isomorphism outside of the exceptional divisor $D=\pi^{-1}(0)$,
$D$ is a normal crossing divisor. Let $E_\sigma$
be irreducible components of the divisor $D$: $D=\bigcup_\sigma E_\sigma$. For a germ
$f\in {\cal O}_{V,0}$, let $v_{\sigma}(f)$ be the multiplicity of the lifting
$\widetilde f=f\circ\pi$ of the germ $f$ to the space $\widetilde V$ of the resolution
along the component $E_\sigma$. The function
$v_\sigma: {\cal O}_{V,0}\setminus\{0\}\to\Z_{\ge0}\cup\{+\infty\}$ defines a valuation (a
{\em divisorial} one) on the field of quotients of the ring ${\cal O}_{V,0}$.

Suppose that the resolution $\pi:(\widetilde V, D)\to (V,0)$ is equivariant with
respect to an action of a finite group $G$ on $(V,0)$, i.e., that the action
of $G$ on $(V,0)$ can be lifted to an action on $(\widetilde V, D)$. In
particular the group $G$ acts on the set $\{E_\sigma\}$ of irreducible
components of the exceptional divisor $D$ by permutations. If the components
$E_{\sigma_1}$ and $E_{\sigma_2}$ are in the same orbit of the action and a
function $f\in {\cal O}_{V,0}$ is equivariant (with respect to a 1-dimensional
representation $\alpha: G\to \C^*$), then $v_{\sigma_1}(f)= v_{\sigma_2}(f)$.

Let us fix $s$ components $E_1$, \dots, $E_s$ of the exceptional divisor $D$.
The corresponding valuations $v_1$, \dots, $v_s$ define a multi-index filtration
on the ring ${\cal O}_{V,0}$ of germs of functions on $(V,0)$ and also a
filtration on the space ${\cal O}^\alpha_{V,0}$ of functions equivariant
with respect to a 1-dimensional representation $\alpha$. Let
$J^\alpha(\vv):=\{f\in {\cal O}^\alpha_{V,0}: v_i(f)\ge v_i, i=1,\ldots,s\}$.
It is easy to see that the filtration $\{J^\alpha(\vv)\}$ is finitely
determined. Let $P^\alpha(\tt)$ ($\tt=(t_1,\ldots,t_s)$) be the Poincar\'e
series of the filtration $\{J^\alpha(\vv)\}$. One has
$$
P^\alpha(\tt)=\int\limits_{\P {\cal O}^\alpha_{V,0}} \tt^{\,\vv(f)}d\chi\,.
$$

\begin{remarks}
{\bf 1.} It is somewhat natural to suppose that for $i\ne j$ the divisors $E_i$ and
$E_j$ are in different orbits of the action of the group $G$, however, from formal point of view this is not important. One can even permit that $E_i=E_j$ for some $i$ and $j$.

\noindent{\bf 2.} In fact it is sufficient to consider the case when $\{E_i\}=\{E_{\sigma}\}$,
i.e. $E_1$, \dots, $E_s$ are all the components of the exceptional divisor $D$.
The Poincar\'e series corresponding to a subset of components of the exceptional
divisor is obtained from the one corresponding to all the components by putting
additional variables equal to $1$.
\end{remarks}

\begin{definition}
The series
$$
P^G(\tt)=\sum\limits_{\alpha\in G^*}\alpha\cdot P^\alpha(\tt)\,\in\,R_1(G)[[t_1,\ldots,t_s]]\,,
$$
where the sum is over all non-isomorphic 1-dimensional representations $\alpha$ of
the group $G$, will be called the {\em equivariant Poincar\'e series} of the divisorial
filtration corresponding to the components $E_1$, \dots, $E_s$ of the exceptional
divisor $D$ of the resolution $\pi:(\widetilde V, D)\to(V, 0)$.
\end{definition}

The discussion above implies that the equivariant Poincar\'e series $P^G(\tt)$ can
be represented as
$$
\int\limits_{\coprod_\alpha \P {\cal O}^\alpha_{V,0}} \alpha(f)\cdot\tt^{\,\vv(f)}d\chi\,,
$$
where $\alpha:\coprod_\alpha \P {\cal O}^\alpha_{V,0}\to R_1(G)$ is the (tautological)
function which sends $\P {\cal O}^\alpha_{V,0}$ to $\alpha\in R_1(G)$.

The ring $R_1(G)$ can be described in the following way. Let the abelian group
$G/(G,G)$ ($(G,G)$ is the commutator of the group $G$) be the direct sum
$\bigoplus\limits_{j=1}^q \Z_{m_j}$ of the cyclic groups $\Z_{m_j}=\Z/m_j\Z$,
let ${\widehat\xi}_j$ be a generator of the group $\Z_{m_j}$, let ${\xi}_j$
be an element of the group $G$ which maps to ${\widehat\xi}_j$ under the
factorization by the commutator $(G,G)$, and let $u_j$ be the
representation of the group $G$ defined by $u_j(\xi_j)=\exp\frac{2\pi i}{m_j}$,
$u_j(\xi_k)=1$ for $k\ne j$. Then the representations
$\prod\limits_{j=1}^q u_j^{\ell_j}$ with $0\le\ell_j<m_j$ form a basis of the ring
$R_1(G)$ as a $\Z$-module and as a ring
$R_1(G)\cong\Z[u_1,\ldots, u_q]/(u_j^{m_j}-1;\ j=1,\ldots, q)$ (see e.g.
\cite{S}).

\section{Divisorial filtrations on the ring ${\cal O}_{\C^2,0}$}\label{sec2}
From now on let $(V,0)=(\C^2,0)$, a finite group $G$ acts on $(\C^2,0)$.
Without loss of generality we suppose this action to be defined by a
representation $G\to \mbox{\bf GL}(2,\C)$.
Let $\pi:(X,D)\to (\C^2,0)$ be a modification of the plane $(\C^2,0)$ by a sequence  of
blow-ups invariant with respect to the group action, that is the action
of the group $G$ on the plane $(\C^2,0)$ lifts to its action on the space
$X$ of the modification. All components $E_\sigma$ of the
exceptional divisor $D$ are isomorphic to the complex projective line
$\CP^1$. Let ${\stackrel{\bullet}{E}}_\sigma$ be the ``smooth part" of the
component $E_\sigma$, i.e. $E_\sigma$ itself minus intersection points with
all other components of the exceptional divisor $D$, let
${\stackrel{\bullet}{D}}=\bigcup\limits_\sigma{\stackrel{\bullet}{E}}_\sigma$
be the smooth part of the exceptional divisor $D$, and let
$\widehat{D}={\stackrel{\bullet}{D}}/G$ be the corresponding factor space,
i.e. the space of orbits of the action of the group $G$ on ${\stackrel{\bullet}{D}}$.

Let $-k_\sigma$ be the self intersection number $E_\sigma\circ E_\sigma$
of the component $E_\sigma$ on the (smooth) surface $X$. Let
$(E_\sigma\circ E_{\sigma^\prime})$ be the intersection matrix of the
components $E_\sigma$ of the exceptional divisor $D$:
$E_\sigma\circ E_\sigma=-k_\sigma$; for $\sigma\ne{\sigma^\prime}$ one has
$E_\sigma\circ E_{\sigma^\prime}=1$ if the components $E_\sigma$ and
$E_{\sigma^\prime}$ intersect (at a point) and $E_\sigma\circ E_{\sigma^\prime}=0$
otherwise. One has $\det(E_\sigma\circ E_{\sigma^\prime})=(-1)^{\#\sigma}$ where
${\#\sigma}$ is the number of components of the exceptional divisor $D$, i.e.
the number of blow-ups which produce the modification $(X,D)$. Let
$M=(m_{\sigma{\sigma^\prime}}):= -(E_\sigma\circ E_{\sigma^\prime})^{-1}$.
All entries of the matrix $M$ are positive integers. The meaning of the entry
$m_{\sigma{\sigma^\prime}}$ is the following. Let ${\widetilde L}_\sigma$
be a germ of a smooth curve on $X$ transversal to the component $E_\sigma$ at a
smooth point of it, i.e. at a point of ${\stackrel{\bullet}{E}}_\sigma$. Let the
curve $L_\sigma=\pi({\widetilde L}_\sigma)\subset(\C^2,0)$ be given by an
equation $h_\sigma=0$, $h_\sigma\in{\cal O}_{\C^2,0}$. Then
$m_{\sigma{\sigma^\prime}}=v_{\sigma^\prime}(h_\sigma)$ ($=v_\sigma(h_{\sigma^\prime})$).
For $\sigma\ne\sigma^\prime$ the entry
$m_{\sigma{\sigma^\prime}}$ is also equal to the intersection number
$L_\sigma\circ L_{\sigma^\prime}$. The diagonal entry $m_{\sigma{\sigma}}$
is equal to $L_\sigma\circ L^\prime _{\sigma }$ where
$L_\sigma=\pi({\widetilde L}_\sigma)$,
$L^\prime_\sigma=\pi({\widetilde L}^\prime_\sigma)$ for smooth curves
${\widetilde L}_\sigma$ and ${\widetilde L}^\prime_\sigma$ transversal
to the component $E_\sigma$ at different smooth points.

For a space $Z$, let $S^kZ=Z^k/S_k$ be the $k$th symmetric power of the space $Z$.
If $Z$ is a quasiprojective variety, $S^kZ$ is also one. For example $S^k\C\cong\C^k$,
$S^k\CP^1\cong\CP^k$.

Let $E_0$ be the first (in the order of blow-ups) component of the exceptional
divisor $D$ ($E_0$ may very well coincide with one of the chosen components
$E_1$, \dots, $E_s$). The component $E_0$ is invariant with respect to the
action of the group $G$. For a point $p$ of $E_0$, i.e. for a line in $\C^2$,
let $\ell_p$ be the corresponding linear function on $\C^2$: the function which
vanishes at $p$; it is well defined up to multiplication by a non-zero constant.
Let $G_p=\{g\in G: gp=p\}$ be the stabilizer of the point $p$ and let $Gp\cong G/G_p$
be the orbit of the point $p$. The homogeneous function
$\prod\limits_{q\in Gp} \ell_q$
of degree equal to the number $\vert Gp\vert = \vert G\vert/\vert G_p\vert$ of points
in the orbit of $p$ is an equivariant one (with respect to a certain 1-dimensional
representation of the group $G$). Let $u_p$ be the corresponding representation.

\begin{examples}
{\bf 1.} Let the group $\Z_m$ act on $(\C^2, 0)$ by $\xi(x,y)=(\xi^k x,
\xi^{\ell} y)$ where $\xi=\exp(2\pi i/m)$ and $\gcd(k, \ell)$ is relatively
prime with $m$ ($k$ and $\ell$ are the weights of the representation of the
group $\Z_m$).
The group $\Z_m^*$ of 1-dimentional representation of the group $\Z_m$ is the
cyclic group of order $m$ generated by the representation $u$ defined by
$u(\xi)=\xi$. The action of the group $\Z_m$ on the projectivization
$\CP^1=(\C^2\setminus\{0\})/\C^*$ of the plane $\C^2$ has two fixed points
$p_1=(1, 0)$ and $p_2=(0, 1)$ and is free outside of them. The linear function
which vanishes at the point $p_1$ is the coordinate $y$. Therefore
$u_{p_1}=u^{-\ell}$. In the same way $u_{p_2}=u^{-k}$. For a point $p\in\CP^1$
different from $p_1$ and $p_2$ one has $u_p=1$.

{\bf 2.} Let $G=D_2^*$ be the subgroup of $\mbox{\bf SL}(\C, 1)$ generated by
the transformations $\sigma(x,y)=(ix, -iy)$ and $\tau(x,y)=(iy,ix)$,
$i=\sqrt{-1}$. The group $G$ is of order
$8$ and has $3$ cyclic subgroups of order $4$ generated by $\sigma$, $\tau$,
and $\tau\sigma$ respectively. Intersection of any two of these subgroups is
the commutator  $(G,G)=\{1, \sigma^2\}$ of the group $G$
($\sigma^2=\tau^2=(\tau\sigma)^2$). The group $G^*$ of 1-dimensional
representations  of the group $G$ is isomorphic to $\Z_2\times\Z_2$ and is
generated by the representations $u_1$ and $u_2$ defined by $u_1(\sigma)=-1$,
$u_1(\tau)=1$,  $u_2(\sigma)=1$, $u_2(\tau)=-1$; $R_1(G)=
\Z[u_1,u_2]/(u_1^2-1, u_2^2-1)$.

The isotropy group of a generic point of
the action of  the group $G$ on the projectivization
$\CP^1=(\C^2\setminus\{0\})/\C^*$ consists of 2 elements and is the commutator
$(G,G)$ of the group $G$. Each cyclic subgroup of order $4$ has two fixed
points which are in the same orbit of the $G$-action. They are $Q_1=(1:0)$ and
$Q_2=(0:1)$ for the subgroup $\langle\sigma\rangle$, $Q_3=(1:1)$ and
$Q_4=(1:-1)$ for the subgroup
$\langle\tau\rangle$, $Q_5=(1:i)$ and $Q_6=(1:-i)$ for the subgroup
$\langle\tau\sigma\rangle$. The pairs $(Q_1, Q_2)$, $(Q_3, Q_4)$, and $(Q_5, Q_6)$ are orbits of the action of the group $G$ on $\CP^1$.

The product of the linear functions vanishing at the points $Q_1$ and $Q_2$ is $xy$. It is invariant with respect to the transformation $\sigma$ and anti-invariant with respect to the transformation $\tau$. Therefore $u_{Q_1}=u_{Q_2}=u_2$. In the same way $u_{Q_3}=u_{Q_4}=u_1$, $u_{Q_5}=u_{Q_6}=u_2$ (the corresponding products of linear functions are $x^2-y^2$ and $x^2+y^2$ respectively). The product of linear functions vanishing at the points of an orbit of the action of the group $G$ different from those consisting of the points $Q_i$ is a homogeneous polynomial of degree $4$ and is invariant with respect to the $G$-action. Therefore for any $Q\ne Q_i$, $i=1,\ldots, 6$, one has $u_Q=1$.
\end{examples}

Let $\{\Xi\}$ be a stratification of the space (in fact of a smooth curve) $\widehat D$
(${\widehat D}=\bigcup\Xi$) such that:
\begin{enumerate}
\item Each stratum $\Xi$ is connected and therefore is
contained in the image of one component $E_\sigma$ under the map
$j:{\stackrel{\bullet}{D}}\to{\widehat D}$ of factorization.
\item Over each strutum $\Xi$ the map $j$ is a covering.
\end{enumerate}
As above, for a point $p\in {\stackrel{\bullet}{D}}$, let $G_p$ be its stabilizer
$\{g\in G: gp=p\}$ and let $Gp$ be the orbit of the point $p$. For all points
from a connected component of the preimage $j^{-1}(\Xi)$ of a stratum $\Xi$
the stabilizer $G_p$ is one and the same and for all points of $j^{-1}(\Xi)$
the stabilizers are conjugate to each other. For $p\in {\stackrel{\bullet}{D}}$,
let $\sigma=\sigma(p)$ be the number of the component $E_\sigma$ which contains
the point $p$. Let
$$
m_{\Xi,i}\ :=\ \sum\limits_{p\in j^{-1}({\widehat p})}
m_{\sigma(p)i}\quad \mbox{for $\widehat p \in \Xi$},
$$
$\mm_{\,\Xi}:=( m_{\Xi,1}, \ldots, m_{\Xi,s})$.

For $p\in{\stackrel{\bullet}{D}}$, let $p_0$ be the corresponding point of the
component $E_0$ of the exceptional divisor: if $p\not\in E_0$, $p_0$ is
the point which belongs to the
closure of the corresponding component of the complement $D\setminus E_0$.
Otherwise $p_0=p$. Let
$$
u_\Xi\ :=\ u_{p_0}^{m_{\sigma(p)0}(\vert Gp\vert/\vert Gp_0\vert)}\,,
$$
where $p\in j^{-1}(\Xi)$. The (1-dimensional) representation $u_{\Xi}$ is a
monomial, say
$\uu^{\,\ll_{\,\Xi}}$, in basic representations,
$\uu=(u_1,\ldots, u_q)$,
$\ll_{\,\Xi}=(\ell_1, \ldots, \ell_q)$.

\begin{theorem}\label{theo1}
\begin{equation}\label{eq2}
{P^G}(\tt)\ =\ \prod\limits_{\Xi}
\left(1-\uu^{\,\ll_{\,\Xi}}\tt^{\,\mm_{\,\Xi}}\right)^{- \chi(\Xi)}\,.
\end{equation}
\end{theorem}

\begin{proof}
Let us fix $\VV\in \Z_{\ge 0}^s$ and suppose that we compute the Poincar\'e
series $P^{G}(\tt)$ up to terms of degree $\VV$ in $\seq{t}1s$. Let
$\P{\cal O}^\alpha_{\C^2,0}(\VV)$ be the set
$\{f\in \P{\cal O}^\alpha_{\C^2,0} : \vv(f)\le \VV \}$. We can make additional
blow-ups of intersection points of components of the exceptional divisor so
that, for any
$f\in\OO_{\C^2,0}$ with
$\vv(f)\le \VV$, the strict transform of the curve $\{f=0\}$ intersects the
exceptional divisor only at smooth points. We shall keep the notation
$(X,D)$ for this modification as well explaining later why new-born components do not
affect the answer.
The space of $G$-invariant effective divisors on ${\stackrel{\bullet}{D}}$, i.e.
of
$G$-invariant
unordered collection of points of ${\stackrel{\bullet}{D}}$ with some multiplicities,
can be in the natural way identified with the space
$$
Y^G\ =\ \prod\limits_{\Xi}\left(\bigcup\limits_{k=0}^\infty S^k \Xi \right) =
\bigcup\limits_{\{k_{\Xi}\}} \left(\prod\limits_{\Xi} S^{k_{\Xi}} \Xi \right)\,.
$$
The space $Y^G$ can be considered as a semigroup with the union of collections
as the semigroup operation.
Let
$$
I_G:\,\coprod\limits_\alpha\,\P{\cal O}^\alpha_{\C^2,0}(\VV)\to Y^G
$$
be the map which sends a (non-zero) function $f\in{\cal O}^\alpha_{\C^2,0}(\VV)$
to
the intersection of the strict transform of the zero-level curve
${\{f=0\}}$ with the exceptional divisor $D$ (a $G$-invariant collection of the intersection points counted with the corresponding multiplicities).

Let $\vv: Y^G \to \Z_{\ge 0}^s$ be the map (a semigroup homomorphism) which
sends a point of the component $\prod\limits_{\Xi} S^{k_{\Xi}} \Xi$ of the space
$Y^G$ to
$\sum\limits_{\Xi} k_{\Xi} \mm_{\,\Xi}$. Let $\uu: Y^G\to G^{*}$ be the map to the
group $G^*$ of one-dimensional representations of the group $G$ which
sends a point of the component $\prod\limits_{\Xi} S^{k_{\Xi}} \Xi$ of the space
$Y^G$ to
$\prod\limits_{\Xi} u_{\Xi}^{k_{\Xi}}$.
One can easily see that $\vv\circ I_G = \vv$ and
$\uu\circ I_G = \uu$ where $\vv$ and $\uu$ in the right hand sides of the
equations are the corresponding maps from
$\coprod\limits_\alpha\,\P{\cal O}^\alpha_{\C^2,0}(\VV)$. The image of the map $I_G$ is the union of all components of the space $Y^G$ with $\vv\le\VV$. For a space (a quasi-algebraic variety) $Z$ one has
$$
\sum\limits_{k=0}^\infty \chi(S^kZ)t^k = (1-t)^{-\chi(Z)}\,.
$$
This implies that
$$
\int\limits_{Y^G} \uu\,\tt^{\,\vv} d\chi = \prod\limits_\Xi \left(1-\uu^{\ll_{\,\Xi}}\tt^{\mm_{\,\Xi}}\right)^{-\chi(\Xi)}\,.
$$
Moreover, all additional strata corresponding to components of the exceptional divisor born under additional blow-ups have the Euler characteristics equal to zero and therefore do not participate in the right hand side of the equation.

For any point $y\in Y^G$ with $\vv(y)\le\VV$, its preimage $I_G^{-1}(y)$ in
$\coprod\limits_\alpha\,\P{\cal O}^\alpha_{\C^2,0}(\VV)$ is an affine space: \cite[Proposition 2]{CDGb}. With the Fubini formula this implies that the integrals of $\uu\,\tt^{\,\vv}$ over the spaces $Y^G$ and $\coprod\limits_\alpha\,\P{\cal O}^\alpha_{\C^2,0}(\VV)$ coincide up to terms of degree $\ge\VV$. Therefore Equation~\ref{eq2} is correct up to terms of degree $\ge\VV$ for any $\VV\in\Z_{\ge 0}^s$. This implies the statement.
\end{proof}

\section{Examples}
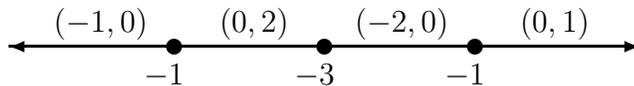
\begin{figure}[h]
$$
\unitlength=1.00mm
\begin{picture}(80.00,20.00)(0,0)
\thicklines
\put(20,10){\line(1,0){40}}
\put(20,10){\circle*{2}}
\put(40,10){\circle*{2}}
\put(60,10){\circle*{2}}
\put(20,10){\vector(-1,0){22}}
\put(60,10){\vector(1,0){22}}
\put(16,5){$-1$}
\put(36,5){$-3$}
\put(56,5){$-1$}
\put(4,12){$(-1,0)$}
\put(26,12){$(0,2)$}
\put(44,12){$(-2,0)$}
\put(66,12){$(0,1)$}
\end{picture}
$$
\caption{Example 1.}
\label{fig1}
\end{figure}

{\bf 1.} Let $G=\Z_3$ act on the plane $(\C^2, 0)$ by $\xi(x,y)=(\xi
x, \xi^{-1} y)$, where $\xi=\exp(2\pi i/3)$ is the generator of the
group $\Z_3$. The factor space $(\C^2, 0)/\Z_3$ is a surface
singularity of type $A_2$. Blowing up the origin we glue-in the
exceptional divisor $\CP^1$ with the self-intersection number
$(-1)$. The action of the group $\Z_3$ on the space of the modification
has two fixed points corresponding to the coordinate axes in
$(\C^2,0)$. At one of them (corresponding to the $x$-axis) the
representation of the group $\Z_3$ has weights $(1,1)$, at the other it has weights
$(-1, -1)$. After blowing up these two points one gets a
modification $X$ of the plane $\C^2$ whose exceptional divisor $D$
consists of three components with the self-intersection numbers
$-1$, $-3$, and $-1$. The information about the space $X$ and the
$\Z_3$-action on it is encoded in Figure~\ref{fig1}.

The vertices correspond to the components of the exceptional
divisor. The arrows correspond to the strict transforms of the
coordinate axes in $(\C^2, 0)$. The numbers under the vertices are
the self-intersection numbers of the corresponding components of the
exceptional divisor. The numbers over each edge ($\mbox{mod } 3$) are
the weights of the representation of the group $\Z_3$ at the
corresponding intersection (fixed) point; each of these numbers
corresponds to the action on the component which it is close to.

The group $\Z_3$ acts trivially on the first and on the third
components of the exceptional divisor (counted from the left). The
action on the second one is non-trivial and has two fixed points.
The fact that at each fixed point of the $\Z_3$-action on the space
$X$ one of the weights of the representation is equal to zero
implies that the factor space $X/\Z_3$ is smooth and is a resolution
of the $A_2$-singularity $(\C^2, 0)/\Z_3$.

\begin{remark}
In general, for an action of a group (say, of a cyclic one) on the
plane $(\C^2, 0)$, it is not possible to get an equivariant
modification $X$ of the plane such that the factor space $X/G$ is
non-singular and therefore is a resolution of the singularity
$(\C^2, 0)/G$. The factor space $X/G$ may have singular points: see
e.g. Example~2 below.
\end{remark}

One has
$$
-\left(E_i\circ E_j\right)= \left(
\begin{array}{ccc}
2 & 1 & 1 \\
1 & 1 & 1 \\
1 & 1 & 2
\end{array}
\right)\,.
$$
As the stratification $\{\Xi\}$ of the factor space
$\widehat{D}={\stackrel{\bullet}{D}}/\Z_3$ one can take the one whose
three strata are isomorphic to the projective line $\CP^1$
minus two points (and therefore have the Euler characteristics equal
to zero) and two strata consists of one point each: the images of
the intersection points of the strict transforms of the coordinate
axes in $\C^2$ with the exceptional divisor $D$. For these two
strata the representations $u_\Xi$ are equal to $u$ and $u^{-1}$
respectively. Therefore
$$
P^{\Z_3}(t_1, t_2, t_3) =
(1-u\,\tt^{(2,1,1)})^{-1}(1-u^{-1}\tt^{(1,1,2)})^{-1}\,.
$$

The exceptional divisor of the resolution
$(X,0)/\Z_3\to(\C^2,0)/\Z_3$ is the image of the exceptional divisor
$D$ of the modification $(X,0)\to(\C^2,0)$ and also consists of three
components. The self-intersection numbers of these components can be
determined in the following way. On the first and on the third
components of the exceptional divisor $D$ the $\Z_3$-action is trivial. Therefore they map
isomorphically on their images. The group $Z_3$ acts only on the
normal bundles to these components. The normal bundles to their
images are the third powers of those over the preimages. Therefore
the self-intersection numbers of the corresponding components of
$D/\Z_3$ are equal to $(-3)$. On the second component of $D$ the
$\Z_3$-action is not trivial. This implies that the
self-intersection number of its image (the factor under the
$\Z_3$-action) is $3$ times less than the self-intersection number
of the preimage and is equal to $(-1)$. If, in $(X,0)/\Z_3$, we
blow-down the second component of the exceptional divisor (with the self-intersection number $(-1)$), we get the standard (minimal) resolution of the
$A_2$-singularity. To get the Poincar\'e series of the divisorial
filtration corresponding to it one should take the part of the
series $P^{\Z_3}(T_1^{1/3}, 1, T_2^{1/3})$ corresponding to the
trivial representation. The reason for the change of variables
$T_1=t_1^3$, $T_2=t_2^3$ is the difference between normal bundles to
the corresponding components of $D$ and to their images described
above. This way one gets
$$
P_{A_2}(T_1,
T_2)=\frac{(1-\TT^{(3,3)})}{(1-\TT^{(2,1)})(1-\TT^{(1,2)})(1-\TT^{(1,1)})}
$$
(cf. \cite{CDGc}).

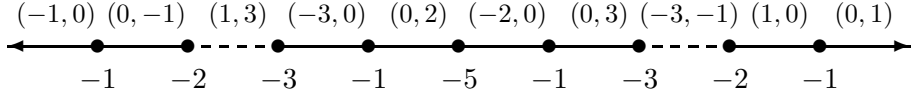
\begin{figure}[h]
$$
\unitlength=0.80mm
\begin{picture}(150.00,20.00)(0,0)
\thicklines
\put(15,10){\line(1,0){15}}
\put(30,10){\dashbox{1.5}(15,0)}
\put(45,10){\line(1,0){60}}
\put(105,10){\dashbox{1.5}(15,0)}
\put(120,10){\line(1,0){15}}

\put(15,10){\circle*{2}}
\put(30,10){\circle*{2}}
\put(45,10){\circle*{2}}
\put(60,10){\circle*{2}}
\put(75,10){\circle*{2}}
\put(90,10){\circle*{2}}
\put(105,10){\circle*{2}}
\put(120,10){\circle*{2}}
\put(135,10){\circle*{2}}

\put(15,10){\vector(-1,0){15}}
\put(135,10){\vector(1,0){15}}

\put(12,3){\small $-1$}
\put(27,3){\small $-2$}
\put(42,3){\small $-3$}
\put(57,3){\small $-1$}
\put(72,3){\small $-5$}
\put(87,3){\small $-1$}
\put(102,3){\small $-3$}
\put(117,3){\small $-2$}
\put(132,3){\small $-1$}

{\footnotesize
\put(0,14){ $(-1,0)$}
\put(15,14){ $(0,-1)$}
\put(32,14){ $(1,3)$}
\put(45,14){ $(-3,0)$}
\put(62,14){ $(0,2)$}
\put(75,14){ $(-2,0)$}
\put(92,14){ $(0,3)$}
\put(103.5,14){ $(-3,-1)$}
\put(122,14){ $(1,0)$}
\put(136,14){ $(0,1)$}
}
\end{picture}
$$
\caption{Example 2.}
\label{fig2}
\end{figure}

\medskip
{\bf 2.} Now let $G=\Z_5$ act on $(\C^2, 0)$ by $\xi(x,y)=(\xi x,
\xi^{-1} y)$ where $\xi=\exp(2\pi i/5)$. The factor space $(\C^2,
0)/\Z_5$ is a surface singularity of type $A_4$. Blowing up the
origin, after that all fixed points of the modification and then all
fixed points of the second one we arrive to a modification the
exceptional divisor of which has 7 irreducible components. The
weights of the $\Z_5$-action at the fixed points which are not
intersection points of the components of the exceptional divisor
(they are intersection points of the exceptional divisor with
the strict transforms of the coordinate axes in $\C^2$) are $(1,1)$ and
$(-1,-1)$ respectively. Blowing up these points one gets the
modification described in Figure~\ref{fig2}.

The factor space of this modification is not smooth: it has two
quotient cyclic singularities of type $(1,3)$ in the terminology
of \cite{P}, corresponding to the intersection points (edges)
depicted by dashed lines. Minimal resolution of such a point has the exceptional
divisor consisting of two components with the self intersection numbers $-2$ and
$-3$ (see \cite{P}). One can easily see that additional
blow-ups at these points do not lead to a modification the factor space of
which is smooth. Just as in Example~1 one has
$$
P^{\Z_5}(t_1, \ldots , t_9) =
(1-u\,\tt^{(4,3,2,3,1,2,1,1,1)})^{-1}(1-u^{-1}\tt^{(1,1,1,2,1,3,2,3,4)})^{-1}\,.
$$
The components of the exceptional divisor $D$ which correspond to the components of
the standard (minimal) resolution of the $A_4$-singularity are
those with numbers 1,4,6, and 9. One can verify this resolving the two singular
points of $X/\Z_5$ and watching the sequence of blow-downs which leads to the
minimal resolution of the $A_4$-singularity. On all of them the
$\Z_5$-action is trivial and therefore they map on the four components of the
resolution of the $A_4$-sigularity isomorphically. The Poincar\'e series of the divisorial
filtration of (the minimal resolution of) the $A_4$-singularity is the part of the
series $P^{\Z_5}(T_1^{1/5}, 1, 1, T_2^{1/5},1, T_3^{1/5}, 1, 1, T_4^{1/5})$ corresponding to the trivial representation. This gives
$$
P_{A_4}(T_1,T_2,T_3,T_4) =
\frac{(1-\TT^{(5,5,5,5)})}{
(1-\TT^{(1,1,1,1)})\,
(1-\TT^{(1,2,3,4)})\,
(1-\TT^{(4,3,2,1)})}\;
$$
(cf \cite{CDGc}).

{\bf 3.} Let $G={\Bbb D}_2^*$ be the subgroup of $\mbox{\bf SL}(\C, 1)$ generated by
the transformations $\sigma(x,y)=(ix, -iy)$ and $\tau(x,y)=(iy,ix)$, $i=\sqrt{-1}$ (see Example~2 of Section~\ref{sec2}). The factor of the plane $\C^2$ by the $G$-action has a
singular point of type $D_4$ (see, e.g., \cite{P}).

\begin{figure}[h]
$$
\unitlength=1.00mm
\begin{picture}(140.00,40)(3,0)
\thicklines
\put(15,18){\line(1,0){15}}
\put(30,18){\line(2,3){8.3}}
\put(30,18){\line(2,-3){8.3}}
\put(30,18){\line(1,0){15}}
\put(30,18){\line(-2,3){8.3}}
\put(30,18){\line(-2,-3){8.3}}

\put(15,18){\vector(-1,0){7.5}}
\put(45,18){\vector(1,0){7.5}}
\put(38.3,30.5){\vector(2,3){4.15}}
\put(21.7,30.5){\vector(-2,3){4.15}}
\put(38.3,5.5){\vector(2,-3){4.15}}
\put(21.7,5.5){\vector(-2,-3){4.15}}

{\thinlines

\put(25,14){\dashbox{1.5}(31,8)}

\put(15,18){\circle*{1}}
\put(22.5,18){\circle*{1}}
\put(30,18){\circle*{1}}
\put(38.3,30.5){\circle*{1}}
\put(34.15,24.25){\circle*{1}}
\put(38.3,5.5){\circle*{1}}
\put(34.15,11.75){\circle*{1}}

\put(45,18){\circle*{1}}
\put(37.5,18){\circle*{1}}
\put(21.7,30.5){\circle*{1}}
\put(25.85,24.25){\circle*{1}}
\put(21.7,5.5){\circle*{1}}
\put(25.85,11.75){\circle*{1}}
}
{\footnotesize
\put(29.5,14.5){$1$}
\put(37,14.5){$2$}
\put(44.5,14.5){$3$}
\put(26.75,24.25){$4$}
\put(22.6,30.25){$5$}
\put(26.75,10.5){$6$}
\put(22.6,4.5){$7$}
}

\thicklines
\put(30,8.5){\vector(1,0){4}}
\put(30,8.5){\vector(-1,0){4}}
\put(20.2,21.4){\vector(2,3){3}}
\put(22.2,24.4){\vector(-2,-3){3}}
\put(39.8,21.4){\vector(-2,3){3}}
\put(37.8,24.4){\vector(2,-3){3}}

\put(58,16.5){\Large $\Rightarrow$}

\put(70,18){\dashbox{2}(20,0)}
\put(90,18){\line(1,0){40}}
\put(70,18){\circle*{2}}
\put(90,18){\circle*{2}}
\put(110,18){\circle*{2}}
\put(110,18){\vector(1,0){20}}

\put(67,13){$-7$}
\put(87,13){$-2$}
\put(107,13){$-1$}

\put(74,20){$(2,-1)$}
\put(95,20){$(1,0)$}
\put(115,20){$(0,1)$}

\put(65,8){\dashbox{1.5}(70,20)}

\end{picture}
$$
\caption{Example 3.}
\label{fig3}
\end{figure}
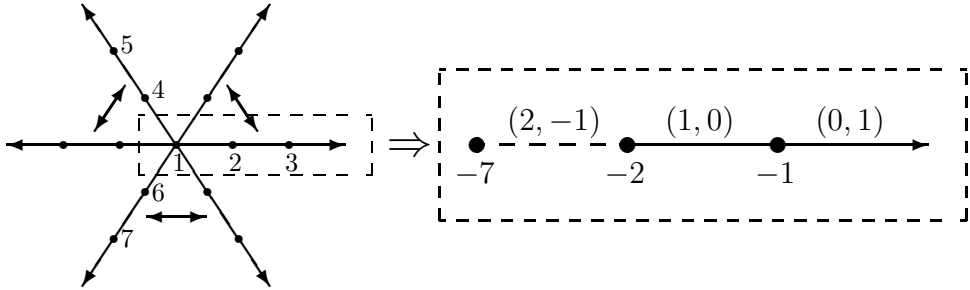

Let us blow-up  the origin in $\C^2$.
On the exceptional divisor (the projective line $(\C^2\setminus\{0\})/\C^*$) there are 6 points $Q_1$, \dots, $Q_6$ whose stabilizers
are cyclic subgroups of $G$ of order $4$ (see Section~\ref{sec2}, Example~2). The stabilizer
of any other point of the exceptional divisor consists of 2 elements and is the commutator
$(G,G)$ of the group $G$. Blowing-up the points $Q_i$, $i=1, \ldots, 6$, one gets $6$ new components of the exceptional divisor. On a component
corresponding to a fixed point $Q_i$ of a cyclic subgroup $H$ of order $4$ the
action of the group $H$ is not trivial. It has two fixed points one of which
is the intersection point with the exceptional divisor of the initial blow-up.
The other one is the intersection point with the strict transform of the line
in $\C^2$ corresponding to the point $Q_i\in(\C^2\setminus\{0\})/\C^*$. Blowing-up these $6$
points which do not
belong to the exceptional divisor of the initial blow-up, one gets $6$ new
components of the exceptional divisor. On the one of them corresponding to a
fixed point of a cyclic subgroup of order $4$ the action of this subgroup is
trivial. We arrive to the modification described in Figure~\ref{fig3}. The
arrows correspond to the strict transforms of the lines in $\C^2$ corresponding
to the points $Q_i$, $i=1, \ldots, 6$. The
numbers in parenthesis are the weights of the action of the corresponding
cyclic group  of order 4 at the corresponding point. Let us choose one
component of the exceptional divisor from each $G$-orbit. (The multiplicities
of any equivariant function on components from the same $G$-orbit are equal.)
These 7 components are indicated on the left hand side of Figure~\ref{fig3}.

A stratification $\{\Xi\}$ of the factor space $\widehat{D}=
{\stackrel{\bullet}{D}}/G$ can be chosen in the following way.
\begin{enumerate}
\item
The factor space $\Xi_1$ of the smooth part of the initial (``central")
component of the exceptional divisor. It has the Euler characteristic equal to
$-1$. The corresponding representation $u_{\Xi_1}$ is trivial.
\item
Three strata consisting of one point each: the images of the intersection points
of the exceptional divisor $D$ with the strict transforms of the lines
corresponding to the points $Q_i$, $i=1,\dots, 6$. The corresponding representations
$u_{\Xi_i}$ are: $u_2$ for the stratum corresponding to the fixed points of the
subgroup $\langle\sigma\rangle$, $u_1$ for that corresponding to the subgroup
$\langle\tau\rangle$, $u_1u_2$ for that corresponding to the subgroup
$\langle\tau\sigma\rangle$.
\item
Six strata with the Euler characteristics equal to zero. They are images of the
smooth parts of 12 components of the exceptional divisor $D$ different from the
central one minus intersection points with the strict transforms of the lines corresponding to the points $Q_i$.
\end{enumerate}

Theorem~\ref{theo1} gives
$$
P^G(\seq t17) = \frac{1- \tt^{(4,4,4,4,4,4,4)}}{
(1-u_2\,\tt^{(2,3,4,2,2,2,2)})
(1-u_1\,\tt^{(2,2,2,3,4,2,2)})
(1-u_1 u_2\,\tt^{(2,2,2,2,2,3,4)})} \, .
$$

The components of the exceptional divisor $D$ of the modification which
corresponds to the components of the standard (minimal) resolution of the $D_4$
singularity $(\C^2,0)/G$ are those with numbers 1, 3, 5, and 7. Therefore
$P_{D_4}(T_1,T_2,T_3,T_4)$ is the part of the series
$P^G(T_1^{1/2},1,T_2^{1/4},1,T_3^{1/4},1,T_4^{1/4})$ corresponding to the
trivial representation. One gets
$$
P_{D_4}(\TT) = \frac{(1-\TT^{(2,1,1,1)})(1+\TT^{(3,2,2,2)})}{
(1-\TT^{(2,2,1,1)})
(1-\TT^{(2,1,2,1)})
(1-\TT^{(2,1,1,2)})
}
$$
(cf \cite{CDGc}).

\section{Filtrations defined by $G$-invariant curves}

Let $(C, 0)\subset(\C^2,0)$ be a plane curve singularity invariant with respect
to the action of the group $G$ on $(\C^2,0)$. Let $C=\bigcup_{i=1}^r C_i$ be the
representation of the curve $C$ as the union of its irreducible components.
Let $\varphi_i:(\C,0)\to(\C^2, 0)$ be the parametrization (uniformization) of
the component $C_i$ of the curve $C$, i.e. $\mbox{Im}\, \varphi_i =C_i$,
$\varphi_i$ is an isomorphism between $\C$ and $C_i$ outside of the origin.
For a germ $f\in \OO_{\C^2,0}$ let $v_i=v_i(f)$ be the power of
the leading term in the power series decomposition of
the germ $f\circ \varphi_i: (\C,0)\to \C$:
$f\circ \varphi_i (\tau) = a \tau^{v_i} + \mbox{terms of higher degree}$,
$a\neq 0$.  If $f\circ \varphi_i (\tau) \equiv 0$, $v_i(f):=\infty$.
The functions (valuations) $v_i$ define a multi-index filtration on
the ring $\OO_{\C^2,0}$ of germs of functions of two variables and
also on the subspace $\OO^\alpha_{\C^2,0}$ of functions equivariant
with respect to a 1-dimensional representation $\alpha$ of the group
$G$. Let $P_C^\alpha(t_1, \ldots, t_r)$ be the Poincar\'e series of
this filtration on the space ${\OO}^\alpha_{\C^2,0}$, let
$$
P_C^G(\tt) :=\sum_\alpha \alpha\cdot P_C^\alpha(\tt)\in
R_1(G)[[t_1, \ldots, t_r]]
$$
(the sum is over all non-isomorphic 1-dimensional representations $\alpha$
of the group $G$) be the {\em equivariant Poincar\'e series} of the filtration
on the ring $\OO_{\C^2,0}$.

\begin{remark}
One can easily see that if the components $C_{i_1}$ and $C_{i_2}$ of the curve
$C$ are in one and the same orbit of the $G$-action on the set of components,
$v_{i_1}(f)=v_{i_2}(f)$ for any $f\in \OO^\alpha_{\C^2,0}$. This means that in all the
monomials of the equivariant Poincar\'e series $ P_C^G(\tt)$ the exponents
at the variables $t_{i_1}$ and $t_{i_2}$ coincide. Therefore without any loss
of generality one can choose one component from each orbit of the $G$-action.
\end{remark}

Let $\pi:(X,D)\to(\C^2,0)$ be a $G$-invariant embedded resolution of the curve
$C$ (i.e. the total transform $\pi^{-1}(C)$ of the curve $C$ is a normal
crossing divisor on the space $X$ of the resolution and the group $G$ acts
on the space (a surface) $X$ as well). Let $D=\bigcup_\sigma E_\sigma$ be the representation
of the exceptional divisor $D=\pi^{-1}(0)$ as the union of its irreducible
components. Let $\{\Xi\}$ be a stratification of the factor space
$\widehat{D}={\stackrel{\bullet}{D}}/G$ described in section~\ref{sec2}.
Let ${\stackrel{\circ}{E}}_\sigma$ be
$E_\sigma$ itself minus intersection points with
all other components of the total transform $\pi^{-1}(C)$ of the curve $C$ (${\stackrel{\circ}{E}}_\sigma\subset{\stackrel{\bullet}{E}}_\sigma$), let
${\stackrel{\circ}{D}}=\bigcup\limits_\sigma{\stackrel{\circ}{E}}_\sigma$,
and let
$\widehat{D}^\prime:={\stackrel{\circ}{D}}/G$ be the space of orbits
of the action of the group $G$ on $\widehat{D}^\prime$.
Let $\{\Xi^\prime\}$ be a stratification of the space $\widehat D^\prime$
with the same properties as the stratification $\{\Xi\}$ of the space
$\widehat D$ in Section~\ref{sec2}.

One has
${\stackrel{\circ}{D}}\subset {\stackrel{\bullet}{D}}$,
$\widehat{D}' \subset \widehat{D}$ and one can suppose that each stratum of the
stratification $\{\Xi'\}$ is a part of a stratum of the stratification
$\{\Xi\}$. For a stratum $\Xi'$, let $\ll_{\,\Xi'}$ be equal to $\ll_{\,\Xi}$ for
the corresponding stratum $\Xi$.
For $i=1, \ldots, r$, let $E_{\sigma(i)}$ be the component of the exceptional
divisor $D$ of the resolution which intersects the component $C_i$ of the curve
$C$.
Let $m_{\Xi',i} := m_{\Xi,\sigma(i)}$,
$\mm_{\,\Xi'} := (m_{\Xi',1}, \ldots, m_{\Xi',r})$.

\begin{theorem}\label{theo2}
$$
P^G_C(\tt)\ =\ \prod\limits_{\Xi^\prime}
\left(1-\uu^{\ll_{\,\Xi^\prime}}\tt^{\mm_{\,\Xi^\prime}}\right)^{-
\chi(\Xi^\prime)}\,.
$$
\end{theorem}

The {\bf proof} is essentially the same as the one of Theorem~\ref{theo1}.

\bigskip

\noindent
Addresses:

\medskip
\noindent
University of
Valladolid
\newline
Dept. of Algebra, Geometry and Topology
\newline
47005
Valladolid, Spain
\newline
E-mail: campillo\symbol{'100}cpd.uva.es,
fdelgado\symbol{'100}agt.uva.es

\medskip
\noindent
Moscow State University
\newline
Faculty of Mathematics and Mechanics
\newline
Moscow, 119992, Russia
\newline
E-mail: sabir\symbol{'100}mccme.ru

\end{document}